\documentclass{article}
\usepackage{amsmath}
\usepackage{amsfonts}
\usepackage{amssymb}
\usepackage{euscript}

\def\dj{d\kern-0.4em\char"16\kern-0.1em}

\newtheorem{thm}{Theorem}
\newtheorem{cor}[thm]{Corollary}

\newtheorem{lem}[thm]{Lemma}
\newtheorem{exmp}{Example}

\date{}

\title{\bf
On $q$--fractional derivatives of Riemann--Liouville and Caputo
type}

\author{
{\bf Miomir S. Stankovi\'c}\\
{\small\it Department of Mathematics, Faculty of Occupational Safety}\\
\small E-mail:\ miomir.stankovic@gmail.com\\
[2mm]
{\bf Predrag M. Rajkovi\'c}\\
{\small\it Department of Mathematics, Faculty of Mechanical Engineering}\\
\small E-mail:\ pedja.rajk@yahoo.com\\
[2mm]
{\bf Sladjana D. Marinkovi\'c}\\
{\small\it Department of Mathematics, Faculty of Electronic Engineering}\\
\small E-mail:\ sladjana@elfak.ni.ac.yu\\
[2mm]
{\bf   University of Ni\v s,\ Serbia}
 }

\begin{document}
\maketitle

{\it Abstract.} Based on the fractional $q$--integral with the
parametric lower limit of integration, we define fractional
$q$--derivative of Riemann--Liouville and Caputo type. The
properties are studied separately as well as relations between
them. Also, we discuss properties of compositions of these
operators.

\smallskip

\smallskip




\section{Introduction}

The fractional differential equations (FDE), as generalizations of
integer-order ones, are used in describing various phenomena in
the science, especially in physics, chemistry and material
science, because of their ability to describe memory effects
\cite{Diethelm}. Today there are a number of concepts with
different definitions of fractional integrals and derivatives and
their applications in various mathematical areas (see, for example
\cite{Podlubny}.

At the first moment, it was considered that it exists unique
definition of fractional derivative until some confusion appeared
in the conclusions. Now, we know that there are two basic types:
Riemann-Liouville and Caputo fractional derivative. Hence two
types of FDE are in use with very important difference in initial
conditions: the first one requires initial conditions for
fractional derivatives; on the contrary, the second one for
integer order derivatives.

Many of continuous scientific problems have their discrete
versions. A way of the treatment is from the point of view of
$q$--calculus (see, for example \cite{Bangerezako}). W.A. Al-Salam
\cite{Al-Salam} and R.P Agarwal \cite{Agarwal} introduced several
types of fractional $q$--integral operators and fractional
$q$--derivatives, always with the lower limit of integration equal
$0$.

However, in some considerations, such as solving of
$q$--differential equation of fractional order with initial values
in nonzero point, it is of interest to allow that the lower limit
of integration is variable. In our paper \cite{MAGT}, we succeed
to generalize this theory in that direction.

In continuation, our purpose in this paper is to define two types
of the fractional $q$--derivatives based on the fractional
$q$--integrals with the parametric lower limit of integration.

\section{Preliminaries}

In the theory of $q$--calculus (see \cite{BHS}), for a real
parameter $q\in \mathbb R^+\setminus\{1\},$ we introduce a
$q$--real number $[a]_q$ by
$$
[a]_q:=\frac{1-q^a}{1-q} \qquad (a \in \mathbb R)\ .
$$
The $q$--analog of the Pochhammer symbol ($q$--shifted factorial)
is defined by:
$$
(a;q)_0=1\ ,\qquad (a;q)_k=\prod_{i=0}^{k-1}(1-aq^i) \qquad
\bigl(k\in\mathbb N\cup \{\infty \}\bigr)\ .
$$
Its natural expansion to the reals  is
\begin{equation}\label{1def}
(a;q)_{\alpha}=\frac{(a;q)_{\infty}}{(aq^{\alpha};q)_{\infty}}
\qquad (\alpha\in \mathbb R)\ .
\end{equation}

Also, $q$--binomial coefficient is given by
\begin{equation}\label{I.42}
{\alpha\brack k}_q=\frac{(q^{-\alpha};q)_k}{(q;q)_k}\ (-1)^k
q^{\alpha k} q^{-\binom k2}\qquad \qquad (k\in\mathbb N,\
\alpha\in\mathbb{R}) \ .
\end{equation}

The following formulas  (see, for example, \cite{BHS},
\cite{Al-SalamVerma} and \cite{MAGT}) will be useful:
\begin{eqnarray}
(a;q)_n&=&\bigl(q^{1-n}/a;q\bigr)_n\ (-1)^n\ a^n\ q^{\binom n2}\ ;
\label{I.7}\\
\frac{(aq^{-n};q)_n}{(bq^{-n};q)_n}&=&\frac{(q/a;q)_n}{(q/b;q)_n}\
\Bigl(\frac ab\Bigr)^n\ ; \label{I.9}\\
(b/a;q)_{\alpha}&=&\sum_{n=0}^{\infty}(-1)^n {{\alpha}\brack{n}}_q
q^{\binom{n}{2}} \Bigl(\frac{b}{a}\Bigr)^n\ ; \label{AlSalam}\\
\frac{(a;q)_{\alpha+n}}{(a;q)_{\alpha}}&=&(aq^\alpha;q)_n \qquad\
\, (n\in\mathbb N;\ a,b,q,\alpha\in\mathbb
R)\label{I11}\ ;\\
\frac{(\mu q^k;q)_{\alpha}}{(\mu;q)_{\alpha}}
&=&\frac{(\mu q^{\alpha};q)_k}{(\mu;q)_k}\qquad\ (\mu,\alpha\in\mathbb R^+)\ ;\label{fracprop2}\\
(q^{k-n};q)_{\alpha}&=&0\qquad\qquad\qquad (k,n\in\mathbb N_0,\
k\le n)\ .\label{fracprop4}
\end{eqnarray}

The next result will have an important role in proving the
semigroup property of the fractional $q$--integral.

Let us denote
\begin{equation}\label{distrib1}
S(\alpha,\beta,\mu)=\sum_{n=0}^{\infty} \frac{(\mu
q^{1-n};q)_{\alpha-1}\ (q^{1+n};q)_{\beta-1}} {(q;q)_{\alpha-1}\
(q;q)_{\beta-1}}\ q^{\alpha n}\ .
\end{equation}
In the paper \cite{MAGT}, the next lemma is proven.
\begin{lem}
For $\mu,\alpha,\beta\in\mathbb R^+$, the following identity is
valid
\begin{equation}\label{distrib}
S(\alpha,\beta,\mu) =\frac{(\mu
q;q)_{\alpha+\beta-1}}{(q;q)_{\alpha+\beta-1}}\ .
\end{equation}
\end{lem}

The $q$--gamma function is defined by
\begin{equation}\label{Gammaq}
\Gamma_q(x)=\frac{(q;q)_\infty}{(q^x;q)_\infty}\ (1-q)^{1-x}
\qquad \bigl(x\in \mathbb R\setminus \{0,-1,-2,\ldots\}\bigr)\ .
\end{equation}
Obviously,
\begin{equation}\label{Gammaq1}
\Gamma_q(x+1)=[x]_q \Gamma_q(x),\qquad \Gamma_q(x)=
(q;q)_{x-1}(1-q)^{1-x}\ .
\end{equation}

The $q$--hypergeometric function \cite{BHS} is defined as
$$
{}_2\phi_1 \Bigl({{a,\ b}\atop{c}}\Big|\ q;x\Bigr)
=\sum_{n=0}^\infty \frac{(a;q)_n(b;q)_n}{(c;q)_n\ (q;q)_n}\ x^n\ .
$$

The $q$--derivative of a function $f(x)$ is defined by
$$
\bigl(D_q f\bigr)(x) = \frac{f(x)-f(q x)}{x-q x}\ \ (x\ne 0)\
,\quad \bigl(D_q f\bigr)(0)=\lim_{x\to 0} \bigl(D_q f\bigr)(x)\ ,
$$
and $q$--derivatives of higher order:
\begin{equation}\label{deriv}
 D^0_q f= f\ , \qquad D^n_q f=D_q\bigl(D^{n-1}_q
f\bigr) \quad (n=1,2,3, \ldots)\ .
\end{equation}

For an arbitrary pair of functions \ $u(x)$  \ and \ $v(x)$ \ and
constants $\alpha, \beta \in\mathbb{R}$, we have linearity and
product rules
$$
\aligned
D_{q}\bigl(\alpha \ u(x)+\beta\ v(x)\bigr)&=\alpha \bigl(D_{q}u\bigr)(x)+\beta \bigl(D_{q}v\bigr)(x),\\
D_{q}\bigl(u(x)\cdot  v(x)\bigr)&=u(qx)\bigl(D_{q}
v\bigr)(x)+v(x)\bigl(D_{q}u\bigr)(x)\ .
\endaligned
$$

In this paper, very  useful examples are the $q$--derivatives of
the next functions:
\begin{eqnarray}
D_q\bigl(x^{\lambda}(a/x;q)_{\lambda}\bigr)&=&[\lambda]_q
x^{\lambda-1}(a/x;q)_{\lambda-1}\ ,\label{xminusa}\\
D_q\bigl(a^{\lambda}(x/a;q)_{\lambda}\bigr)&=&-[\lambda]_q
a^{\lambda-1}(qx/a;q)_{\lambda-1}\ ,\label{aminusx}\\
D_q\bigl(x^{\lambda}\bigr)&=&[\lambda]_q x^{\lambda-1}\
.\label{x^lambda}
\end{eqnarray}

The $q$-integral is defined by
\begin{equation} \label{(qint0)}
\bigl(I_{q,0}f\bigr)(x)=\int_0^x f(t)\,d_q t
=x(1-q)\sum_{k=0}^\infty f(xq^k)\,q^k \quad (0\le |q|< 1),
\end{equation}
and
\begin{equation} \label{(qint)}
\bigl(I_{q,a}f\bigr)(x)=\int_a^x f(t)\,d_q t =\int_0^x f(t)\,d_q
t-\int_0^a f(t)\,d_q t.
\end{equation}
However, these definitions cause troubles in research as they
include the  points outside of the interval of integration (see
\cite{Gauchman}). In the case when the lower limit of integration
is $a=xq^n$, i.e., when it is determined for some choice of  $x$,
$q$ and positive integer $n$, the $q$--integral (\ref{(qint)})
becomes
\begin{equation}
\int_{xq^n}^x f(t)\ d_q t= x(1-q)\sum_{k=0}^{n-1}f(xq^k)q^k \ .
\label{(rest_qint)}
\end{equation}

As for $q$-derivative, we can define $I_{q,a}^n$ operator  by
$$
I_{q,a}^0f=f,\qquad
I_{q,a}^nf=I_{q,a}\bigl(I_{q,a}^{n-1}f\bigr)\quad
(n=1,2,3,\ldots)\ .
$$
For $q$--integral and $q$--derivative operators the following is
valid:
$$
\bigl(D_q I_{q,a} f\bigr)(x)=f(x)\ ,\qquad \bigl(I_{q,a}D_q f\bigr)(x)=f(x)-f(a)\ ,
$$
and, more generally,
\begin{equation}\label{DerInt}
\bigl(D_q^n I_{q,a}^n f\bigr)(x)=f(x)\ \quad (n\in\mathbb N)\ ,
\end{equation}
\begin{equation}\label{IntDer}
\bigl(I_{q,a}^n D_q^n f\bigr)(x)=f(x)
-\sum_{k=0}^{n-1}\frac{\bigl(D_q^k f\bigr)(a)}{[k]_q!}\ x^k(a/x;q)_k
\quad (n\in\mathbb N)\ .
\end{equation}

The formula for $q$--integration by parts is
\begin{equation}\label{IntParts}
\int_a^b u(x)\bigl(D_q v\bigr)(x)\,d_q x
=\bigl[u(x)v(x)\bigr]_a^b-\int_a^b v(qx)\bigl(D_q u\bigr)(x)\,d_q
x\ .
\end{equation}

\section{ The fractional $q$--integral}

In all further considerations we assume that the functions are
defined in  an interval $(0,b)$\ $(b>0)$, and $a\in(0,b)$ is  an
arbitrary fixed point. Also, the required $q$--derivatives and
$q$--integrals exist and the convergence of the series mentioned
in the proofs is assumed.

\smallskip
\noindent {\bf Definition 1}\ The  {\it fractional $q$--integral }
is
\begin{equation}\label{fracAgain}
\bigl(I_{q,a}^{\alpha} f\bigr)(x)
=\frac{x^{\alpha-1}}{\Gamma_q(\alpha)}\int_a^x
(qt/x;q)_{\alpha-1}\ f(t)\ d_q t \qquad (a<x;\ \alpha\in \mathbb
R^{+})\ .
\end{equation}

\begin{lem}
The fractional $q$--integral $(\ref{fracAgain})$ can be written in
the equivalent form
\begin{equation}\label{fracAgainSt}
\bigl(I_{q,a}^{\alpha} f\bigr)(x) = \int_a^x \ f(t)\ d_{q}
w_\alpha(x,t) \qquad (\alpha\in \mathbb R^{+})\ ,
\end{equation}
where $w_\alpha(x,t)$ is the function defined by
\begin{equation}
\label{Weight} w_\alpha(x,t)=\frac{1}{\Gamma_q(\alpha+1)}
 \bigl(x^\alpha - x^\alpha(t/x;q)_{\alpha}\bigr) \qquad (\alpha\in \mathbb R^{+})\ .
\end{equation}
\end{lem}
\noindent{\it Proof.} It is enough to notice that the
$q$--differential of $w_\alpha(x,t)$ over variable $t$ is
\begin{equation}
d_{q} w_\alpha(x,t) = D_{q}\ w_\alpha(x,t)\ d_{q} t
=\frac{x^{\alpha-1}(qt/x;q)_{\alpha-1}}{\Gamma_q(\alpha)}\ d_{q} t
\ . \Box
\end{equation}
Using formula (\ref{AlSalam}), the integral (\ref{fracAgain}) can
be written as
\begin{equation}\label{frac1}
\bigl(I_{q,a}^{\alpha} f\bigr)(x)
=\frac{x^{\alpha-1}}{\Gamma_q(\alpha)} \sum_{k=0}^{\infty} (-1)^k
{{\alpha-1}\brack {k}}_q q^{\binom{k+1}{2}}x^{-k}\int_a^x t^k
f(t)\ d_q t \quad (\alpha\in \mathbb R^{+})\ .
\end{equation}

Putting $\alpha=1$ in (\ref{frac1}), we get $q$--integral
(\ref{(qint)}).

The fractional integral (see, for example \cite{Podlubny}) is the
limitary case of (\ref{fracAgain}) when $q$ arises to  $1$,  since
$$
\lim_{q\nearrow 1}x^{\alpha-1}(qt/x;q)_{\alpha-1}=(x-t)^{\alpha-1}
\ .
$$

 Obviously, the next equality holds:
\begin{equation}\label{qFInt0}
\bigl(I_{q,a}^{\alpha}f\bigr)(a)=
\frac{a^{\alpha-1}}{\Gamma_q(\alpha)}\int_a^a
(qt/a;q)_{\alpha-1}f(t)\ d_q t=0 \ .
\end{equation}

\begin{lem}\label{Ialbasic}
For $\alpha\in \mathbb R^{+}$, the following is valid:
$$
\bigl(I_{q,a}^{\alpha}f\bigr)(x)=\bigl(I_{q,a}^{\alpha+1}D_q
f\bigr)(x)
+\frac{f(a)}{\Gamma_q(\alpha+1)}x^{\alpha}(a/x;q)_{\alpha}\qquad
(a<x)\ .
$$
\end{lem}

\noindent{\it Proof.} According to the formula (\ref{aminusx}),
the $q$--derivative over the variable $t$ is
$$
D_{q}\bigl(x^{\alpha}(t/x;q)_{\alpha}\bigr)=-[\alpha]_q
x^{\alpha-1}(qt/x;q)_{\alpha-1}\ .
$$
Using the $q$--integration by parts (\ref{IntParts}), we obtain
$$
\aligned \bigl(I_{q,a}^{\alpha} f\bigr)(x)
&=-\frac{1}{[\alpha]_q\Gamma_q(\alpha)}\int_a^x
D_q\bigl(x^{\alpha}(t/x;q)_{\alpha}\bigr)
f(t)d_q t\\
&=\frac{1}{\Gamma_q(\alpha+1)}\Bigl(x^{\alpha}(a/x;q)_{\alpha}f(a)+
\int_a^x x^{\alpha}(qt/x;q)_{\alpha}\bigl(D_q f\bigr)(t) d_q t\Bigr)\\
&=\bigl(I_{q,a}^{\alpha+1}D_q
f\bigr)(x)+\frac{f(a)}{\Gamma_q(\alpha+1)}x^{\alpha}(a/x;q)_{\alpha}\
. \quad \Box
\endaligned
$$

\begin{lem} For $\alpha,\beta\in \mathbb R^{+}$, the following is valid:
$$
\int_0^a (qt/x;q)_{\beta-1}\bigl(I_{q,a}^{\alpha}f\bigr)(t) d_q
t=0\qquad (a<x) \ .
$$
\end{lem}

\noindent{\it Proof.} Using formulas (\ref{fracprop4}) and
(\ref{(rest_qint)}), for $n\in\mathbb N_0$, we have
$$
\aligned
\bigl(I_{q,a}^{\alpha}f\bigr)(aq^n)
&=\frac{1}{\Gamma_q(\alpha)}\int_a^{aq^n} (aq^n)^{\alpha-1}
\bigl((qu)/(aq^n);q\bigr)_{\alpha-1} f(u)d_q u\\
&=\frac{-a^\alpha(1-q)}{\Gamma_q(\alpha)}\sum_{j=0}^{n-1}(q^n)^{\alpha-1}
(q^{j+1-n};q)_{\alpha-1} f(aq^j)q^j =0\ .
\endaligned
$$
From the other side, according to the definition of $q$--integral, we have
$$
\int_0^a  (qt/x;q)_{\beta-1}
\bigl(I_{q,a}^{\alpha}f\bigr)(t) d_q t
=a(1-q)\sum_{n=0}^{\infty}(aq^{n+1}/x;q)_{\beta-1}\bigl(I_{q,a}^{\alpha}f\bigr)(aq^n)q^n\ ,
$$
what is obviously equal to zero\ \ . $\Box$

\begin{thm} \label{kom}
Let $\alpha,\beta\in\mathbb R^+$. The $q$--fractional integration
has the following semigroup pro\-per\-ty
$$
\bigl(I_{q,a}^{\beta}I_{q,a}^{\alpha}f\bigr)(x)
=\bigl(I_{q,a}^{\alpha+\beta}f\bigr)(x)\qquad (a<x)\ .
$$
\end{thm}

\noindent{\it Proof.} By previous lemma, we have
$$
\bigl(I_{q,a}^{\beta}I_{q,a}^{\alpha}f\bigr)(x)
=\frac{x^{\beta-1}}{\Gamma_q(\beta)}\int_0^x
(qt/x;q)_{\beta-1}\bigl(I_{q,a}^{\alpha}f\bigr)(t) d_q t,
$$
i.e.,
$$
\aligned
\bigl(I_{q,a}^{\beta}I_{q,a}^{\alpha}f\bigr)(x)
 &=\frac{x^{\beta-1}}{\Gamma_q(\alpha)\Gamma_q(\beta)}\int_0^x
(qt/x;q)_{\beta-1}\ t^{\alpha-1}\int_0^t (qu/t;q)_{\alpha-1}\ f(u) d_q u \\
&- \frac{x^{\beta-1}}{\Gamma_q(\alpha)\Gamma_q(\beta)}\int_0^x
(qt/x;q)_{\beta-1}\ t^{\alpha-1} \int_0^a (qu/t;q)_{\alpha-1}\
f(u) d_q u\ .
\endaligned
$$
Since, as it was proven in the paper \cite{Agarwal}, the equality
$$
\bigl(I_{q,0}^{\beta}I_{q,0}^{\alpha}f\bigr)(x) = \bigl(I_{q,0}^{\alpha+\beta}f\bigr)(x)
$$
is valid, we conclude that
$$
\aligned
 \bigl(I_{q,a}^{\beta}I_{q,a}^{\alpha}f\bigr)(x)
&=\bigl(I_{q,0}^{\alpha+\beta}f\bigr)(x)\\
&- \frac{x^{\beta-1}}{\Gamma_q(\alpha)\Gamma_q(\beta)}\int_0^x
(qt/x;q)_{\beta-1}\ t^{\alpha-1} \int_0^a (qu/t;q)_{\alpha-1}\
f(u) d_q u\ .
\endaligned
$$
Furthermore, we can write
$$
\aligned \big(I_{q,a}^{\beta}I_{q,a}^{\alpha}f\big)(x)
&=\big(I_{q,a}^{\alpha+\beta}f\big)(x)
+\frac{x^{\alpha+\beta-1}}{\Gamma_q(\alpha+\beta)}\int_0^a
(qt/x;q)_{\alpha+\beta-1}f(t) d_q t\\
 &- \frac{x^{\beta-1}}{\Gamma_q(\alpha)\Gamma_q(\beta)}\int_0^x
(qt/x;q)_{\beta-1} t^{\alpha-1} \int_0^a (qu/t;q)_{\alpha-1}f(u)
d_q u,
\endaligned
$$
wherefrom
$$
\big( I_{q,a}^{\beta}I_{q,a}^{\alpha}f\big)(x)
=\big(I_{q,a}^{\alpha+\beta}f\big)(x) +a(1-q)\sum_{j=0}^{\infty} c_j f(aq^j)
q^j,
$$
with
$$
\aligned c_j&=
\frac{x^{\alpha+\beta-1}(aq^{j+1}/x;q)_{\alpha+\beta-1}}{\Gamma_q(\alpha+\beta)}\\
&\qquad
-\frac{x^{\alpha+\beta-1}(1-q)}{\Gamma_q(\alpha)\Gamma_q(\beta)}\sum_{n=0}^{\infty}
(q^{n+1};q)_{\beta-1}\ q^{n(\alpha-1)}(aq^{j+1-n}/x;q)_{\alpha-1}\
q^n.
\endaligned
$$
By using  formulas (\ref{fracprop2}) and (\ref{Gammaq}), we get
$$
\aligned
c_j&=\bigl((1-q)x\bigr)^{\alpha+\beta-1}\\
&\qquad \times \left\{
\frac{(aq^{j+1}/x;q)_{\alpha+\beta-1}}{(q;q)_{\alpha+\beta-1}}
-\sum_{n=0}^{\infty}
\frac{(q^{n+1};q)_{\beta-1}}{(q;q)_{\beta-1}}\
\frac{(aq^{j+1-n}/x;q)_{\alpha-1}}{(q;q)_{\alpha-1}}\
q^{n\alpha}\right\}\ .
\endaligned
$$
Putting $\mu=q^ja/x$ into (\ref{distrib}), we see that $c_j=0$ for
all $j\in \mathbb N$, which completes the proof. $\Box$

\begin{cor}\label{DIal}
For $\alpha\ge n\ (n\in\mathbb N)$ the following is valid:
$$
\bigl(D_q^nI_{q,a}^{\alpha}f\bigr)(x)=\bigl(I_{q,a}^{\alpha-n}f\bigr)(x)\qquad
(a<x) \ .
$$
\end{cor}

\noindent{\it Proof.} The statement follows from Theorem~\ref{kom}
and property (\ref{DerInt}).\ \ $\Box$

\section{The fractional $q$--derivative of \\Riemann-Liouville
type}

On the basis  of fractional $q$--integral, we can define
$q$--derivative of real order.

\smallskip

\noindent {\bf Definition 2}\ \ The {\it fractional
$q$--derivative of Riemann--Liouville type} is
\begin{equation}\label{Rie-Liou}
\bigl(D_{q,a}^{\alpha} f\bigr)(x) =\left\{
\begin{array}{cc}
\bigl(I_{q,a}^{-\alpha}f\bigr)(x), & \alpha \le 0\\ \\
\Bigl(D^{\lceil\alpha\rceil}_q
I_{q,a}^{\lceil\alpha\rceil-\alpha}f\Bigr)(x), &  \alpha> 0,
\end{array}\right.
\end{equation}
where $\lceil\alpha\rceil$ denotes the smallest integer greater or
equal to $\alpha$.

\smallskip

Notice that $\bigl(D_{q,a}^{\alpha} f\bigr)(x)$ has subscript $a$
to emphasize that it depends on the lower limit of integration
used in definition (\ref{Rie-Liou}). Since $\lceil\alpha\rceil$ is
a positive integer for $\alpha\in\mathbb R^+$, then for
$\bigl(D^{\lceil\alpha\rceil}_q f\bigr)(x)$ we apply definition
(\ref{deriv}).

According to definition and (\ref{qFInt0}), we can easily prove
that
\begin{equation}\label{Dqf(a)}
\bigl(D_{q,a}^{\alpha} f\bigr)(a)=0\qquad (\forall
\alpha\in\mathbb R\setminus\mathbb N_0)\ .
\end{equation}

\begin{thm}\label{DRalmE}
 For $\alpha\in \mathbb R$, the following is
 valid:
$$
\bigl(D_{q} D_{q,a}^{\alpha} f\bigr)(x)=\bigl(D_{q,a}^{\alpha+1}
f\bigr)(x)\qquad (a<x)\ .
$$
\end{thm}

\noindent{\it Proof.} According to the formula (\ref{deriv}), the
statement is true for $\alpha\in\mathbb N_0$. For others, we will
consider three cases.

For $\alpha \le -1$, according to Theorem \ref{kom}, we have
$$
\aligned
\bigl(D_{q} D_{q,a}^{\alpha} f\bigr)(x)&=\bigl(D_{q}I_{q,a}^{-\alpha} f\bigr)(x)
=\bigl(D_{q} I_{q,a}^{1-\alpha-1} f\bigr)(x)\\
&=\bigl(D_{q}I_{q,a}I_{q,a}^{-\alpha-1}f\bigr)(x)
=\bigl(I_{q,a}^{-(\alpha+1)}f\bigr)(x) =\bigl(D_{q,a}^{\alpha+1} f\bigr)(x).
\endaligned
$$
In the case $-1< \alpha <0$, i.e.,  $0< \alpha +1 <1$, we obtain
$$
\bigl(D_{q}D_{q,a}^{\alpha} f\bigr)(x)=\bigl(D_{q}I_{q,a}^{-\alpha} f\bigr)(x)
=\bigl(D_{q}I_{q,a}^{1-(\alpha+1)} f\bigr)(x)=\bigl(D_{q,a}^{\alpha+1} f\bigr)(x).
$$
At last, if $\alpha=n+\varepsilon$, $n\in\mathbb N_0$, $0< \varepsilon<1$,
then $\alpha+1\in(n+1,n+2)$, so we get
$$
(D_{q}D_{q,a}^{\alpha} f)(x) =(D_{q}D^{n+1}_q
I_{q,a}^{1-\varepsilon}f)(x)
=(D^{n+2}_q I_{q,a}^{1-\varepsilon}f)(x)
=(D_{q,a}^{\alpha+1}f)(x).\ \Box
$$

\begin{thm}\label{DRalmD}
 For $\alpha\in \mathbb R\setminus\mathbb N_0$, the following is
 valid:
$$
\bigl(D_{q,a}^{\alpha}D_q f\bigr)(x)
=\bigl(D_{q,a}^{\alpha+1}f\bigr)(x)-\frac{f(a)}{\Gamma_q(-\alpha)}x^{-\alpha-1}(a/x;q)_{-\alpha-1}\quad
(a<x)\ .
$$
\end{thm}

\noindent{\it Proof.}
Let us consider two cases. If $\alpha < 0$,
then, with respect to Lemma \ref{Ialbasic}, Theorem~\ref{kom} and
formulas (\ref{xminusa}) and (\ref{DerInt}), the following holds:
$$
\aligned
\bigl(D_{q,a}^{\alpha+1}f\bigr)(x)&=\bigl(D_{q}D_{q,a}^{\alpha}f\bigr)(x)
=\bigl(D_{q}I_{q,a}^{-\alpha}f\bigr)(x)\\
&=D_{q}\Bigl(\bigl(I_{q,a}^{-\alpha+1}D_qf\bigr)(x)+
\frac{f(a)}{\Gamma_q(-\alpha+1)}\ x^{-\alpha}(a/x;q)_{-\alpha}\Bigr)\\
&=\bigl(D_{q}I_{q,a}I_{q,a}^{-\alpha}D_qf\bigr)(x)
+\frac{f(a)}{\Gamma_q(-\alpha+1)}\ [-\alpha]_q x^{-\alpha-1}(a/x;q)_{-\alpha-1}\\
&=\bigl(D_{q,a}^{\alpha}D_qf\bigr)(x)+\frac{f(a)}{\Gamma_q(-\alpha)}x^{-\alpha-1}(a/x;q)_{-\alpha-1}.
\endaligned
$$
If $\alpha>0$, there exist $n\in\mathbb N_0$ and $\varepsilon\in
(0,1)$, such that $\alpha = n+\varepsilon$. Then, applying the
similar procedure, we get
$$
\aligned
\bigl(D_{q,a}^{\alpha+1}f\bigr)(x)&=\bigl(D_{q}D_{q,a}^{\alpha}f\bigr)(x)
=\bigl(D_{q}D_q^{n+1}I_{q,a}^{1-\varepsilon}f\bigr)(x)\\
&=D_{q}^{n+2}\Bigl(\bigl(I_{q,a}^{2-\varepsilon}D_qf\bigr)(x)
+\frac{f(a)}{\Gamma_q(2-\varepsilon)}
\ x^{1-\varepsilon}(a/x;q)_{1-\varepsilon}\Bigr)\\
&=\bigl(D_{q}^{n+1}D_qI_{q,a}I_{q,a}^{1-\varepsilon}D_qf\bigr)(x)
+\frac{f(a)}{\Gamma_q(2-\varepsilon)}
D_{q}^{n+2}\bigl(x^{1-\varepsilon}(a/x;q)_{1-\varepsilon}\bigr)\\
&=\bigl(D_{q}^{n+1}I_{q,a}^{1-\varepsilon}D_qf\bigr)(x)
+\frac{f(a)}{\Gamma_q(-\varepsilon-n)}
\bigl(x^{-\varepsilon-n-1}(a/x;q)_{-\varepsilon-n-1}\bigr)\\
&=\bigl(D_{q,a}^{\alpha}D_qf\bigr)(x)+\frac{f(a)}{\Gamma_q(-\alpha)}
x^{-\alpha-1}(a/x;q)_{-\alpha-1}.\ \Box
\endaligned
$$

\begin{cor}
The semigroup property for fractional $q$--derivative of
Riemann--Liouville type is not valid, i.e., in general
$$
\bigl(D_{q,a}^{\alpha} D_{q,a}^{\beta} f\bigr)(x) \ne
\bigl(D_{q,a}^{\alpha+\beta}f\bigr)(x)\ .
$$
\end{cor}

\smallskip

\begin{exmp}\rm
Notice that from
$$
D_{q,a}^{n+\varepsilon}
\bigl(x^{\varepsilon-1}(a/x;q)_{\varepsilon-1}\bigr)= 0 \qquad
(n\in\mathbb N_0;\ 0<\varepsilon<1)
$$
we have two different conclusions. From one side, it is true
$$
\lim_{\varepsilon \to 1}D_{q,a}^{n+\varepsilon}
\bigl(x^{\varepsilon-1}(a/x;q)_{\varepsilon-1}\bigr)=0=
\bigl(D_{q}^{n+1}{\bf 1}\bigr)(x)=
D_{q}^{n+1}\bigl(x^{0}(a/x;q)_{0}\bigr) \ .
$$
But, from the other side, it is
$$
\lim_{\varepsilon \to 0}D_{q,a}^{n+\varepsilon}
\bigl(x^{\varepsilon-1}(a/x;q)_{\varepsilon-1}\bigr) =0\ne
D_{q}^{n}\bigl(x^{-1}(a/x;q)_{-1}\bigr)\ .
$$
So, we conclude that the mapping $\alpha\mapsto D_{q,a}^{\alpha}
f$ is not continuous from the right side over variable $\alpha$.

\end{exmp}

\section{The fractional $q$--derivative of Caputo type}

If we change the order of operators, we can introduce another type
of fractional $q$--derivative.

\smallskip

\noindent {\bf Definition 2}\ \ The {\it the fractional
$q$-derivative of Caputo type} is
\begin{equation}\label{Caputo}
\bigl({}_{\star}D_{q,a}^{\alpha} f\bigr)(x) =\left\{
\begin{array}{cc}
\bigl(I_{q,a}^{-\alpha}f\bigr)(x), & \alpha \le 0\\ \\
\Bigl(I_{q,a}^{\lceil\alpha\rceil-\alpha}D^{\lceil\alpha\rceil}_q
f(x)\Bigr), &  \alpha> 0.
\end{array}\right.
\end{equation}

\begin{thm}\label{DCalmE}
 For $\alpha\in \mathbb R\setminus\mathbb N_0$ and $a<x$, the following is
 valid:
$$
\bigl({}_{\star}D_{q,a}^{\alpha+1} f\bigr)(x) -\bigl(
{}_{\star}D_{q,a}^{\alpha} D_{q}f\bigr)(x)= \left\{
\begin{array}{cc}
\dfrac{f(a)}{\Gamma_q(-\alpha)}\ x^{-\alpha-1}(a/x;q)_{-\alpha-1}\ , & \alpha\le -1\ ,\\
\\
0\ , & \alpha>-1\ .
\end{array}
\right.
$$
\end{thm}

\noindent{\it Proof.} As in Theorem \ref{DRalmE}, we will consider
three cases. For $\alpha < -1$, according to Lemma \ref{Ialbasic},
we have
$$
\aligned \bigl({}_{\star}D_{q,a}^{\alpha+1}
f\bigr)(x)&=\bigl(I_{q,a}^{-\alpha-1}f\bigr)(x)
=\bigl(I_{q,a}^{-\alpha}D_{q} f\bigr)(x)
+\dfrac{f(a)}{\Gamma_q(-\alpha)}\
x^{-\alpha-1}(a/x;q)_{-\alpha-1}\\
&=\bigl({}_{\star}D_{q,a}^{\alpha}D_{q}f\bigr)(x)
+\dfrac{f(a)}{\Gamma_q(-\alpha)}\
x^{-\alpha-1}(a/x;q)_{-\alpha-1}\ .
\endaligned
$$
In the case $-1< \alpha <0$, i.e.,  $0< \alpha +1 <1$, we obtain
$$
\bigl({}_{\star}D_{q,a}^{\alpha+1} f\bigr)(x)
=\bigl(I_{q,a}^{1-(\alpha+1)}D_{q} f\bigr)(x)
=\bigl(I_{q,a}^{-\alpha}D_{q} f\bigr)(x)
=\bigl({}_{\star}D_{q,a}^{\alpha}D_{q}f\bigr)(x)
$$
Finally, if $\alpha=n+\varepsilon$, $n\in\mathbb N_0$, $0<
\varepsilon<1$, then $\alpha+1\in(n+1,n+2)$, so we get
$$
\aligned ({}_{\star}D_{q,a}^{\alpha+1}f)(x)
&=\bigl(I_{q,a}^{1-\varepsilon}D^{n+2}_q f\bigr)(x)\\
&=\bigl(I_{q,a}^{1-\varepsilon}D^{n+1}_q D_q f\bigr)(x)
=({}_{\star}D_{q,a}^{\alpha}D_q f)(x)\ . \ \Box
\endaligned
$$

\begin{thm}\label{DCalmD}
 For $\alpha\in \mathbb R\setminus\mathbb N_0$ and $a<x$, the following is
 valid:
$$
\aligned \bigl(D_{q}\ {}_{\star}D_{q,a}^{\alpha}f\bigr)(x)& -
\bigl({}_{\star}D_{q,a}^{\alpha+1} f\bigr)(x)\\
&= \left\{
\begin{array}{cc}
0\ , & \alpha<-1\ ,\\    \\
\dfrac{\bigl(D_{q}^{\lceil\alpha\rceil}f\bigr)(a)}{\Gamma_q(\lceil\alpha\rceil-\alpha)}
\
x^{\lceil\alpha\rceil-\alpha-1}(a/x;q)_{\lceil\alpha\rceil-\alpha-1}\
, & \alpha>-1\ .
\end{array}
\right.
\endaligned
$$
\end{thm}

\noindent{\it Proof.}  At first, let $\alpha < 0$. Using Lemma
\ref{Ialbasic}, Theorem \ref{kom} and formulas (\ref{xminusa}) and
(\ref{DerInt}), we get
$$
\aligned \bigl(D_{q}\ {}_{\star}D_{q,a}^{\alpha}f\bigr)(x)
&=\bigl(D_{q}\ I_{q,a}^{-\alpha}f\bigr)(x)\\
&=\bigl(D_{q}I_{q,a}^{-\alpha+1}D_qf\bigr)(x)+\frac{f(a)}{\Gamma_q(-\alpha+1)}
D_{q}\Bigl(x^{-\alpha}(a/x;q)_{-\alpha}\Bigr)\\
&=\bigl({}_{\star}D_{q,a}^{\alpha}D_qf\bigr)(x)
+\frac{f(a)}{\Gamma_q(-\alpha)}\ x^{-\alpha-1}(a/x;q)_{-\alpha-1}\\
\endaligned
$$
The required equalities are valid both for $\alpha<-1$ or
$-1<\alpha<0$, according to Lemma \ref{DCalmE}.

If $\alpha >0$, there exist $n\in\mathbb N_0$ and $\varepsilon\in
(0,1)$, such that $\alpha = n+\varepsilon$. Then, applying the
similar procedure, we get
$$
\aligned \bigl(D_{q}\ {}_{\star}D_{q,a}^{\alpha}f\bigr)(x)
&=\bigl(D_{q}I_{q,a}^{1-\varepsilon}D_q^{n+1}f\bigr)(x)\\
&=\bigl(D_{q}I_{q,a}^{2-\varepsilon}D_q^{n+2}f\bigr)(x)
+\frac{\bigl(D_q^{n+1}f\bigr)(a)}{\Gamma_q(2-\varepsilon)}
D_{q}\Bigl(x^{1-\varepsilon}(a/x;q)_{1-\varepsilon}\Bigr)\\
&=\bigl({}_{\star}D_{q,a}^{\alpha+1}f\bigr)(x)
+\frac{D_q^{n+1}f(a)}{\Gamma_q(n+1-\alpha)}x^{n-\alpha}(a/x;q)_{n-\alpha}\
. \ \Box
\endaligned
$$


\section{The fractional $q$--integrals and $q$--derivatives \\of some elementary functions}

We will use previous results to evaluate fractional $q$--integrals
and $q$--derivatives of some well-known functions in explicit
form. Here, it is very useful to remind on $q$--form of Taylor theorem
\begin{equation}\label{qTaylor}
f(x)=\sum_{k=0}^{\infty}\frac{(D_q^k
f)(a)}{[k]_q!}x^{k}(a/x;q)_{k}
\end{equation}
given by Jackson (see \cite{Al-SalamVerma}). The next lemma will have crucial role in
reaching of our goal.

\begin{lem}\label{Irec}
For $\alpha\in\mathbb R^{+}\setminus \mathbb N_0$,\ $\lambda\in
(-1,\infty)$, the following is valid
\begin{eqnarray*}
I_{q,a}^{\alpha}\bigl(x^{\lambda}(a/x;q)_{\lambda}\bigr)
&=&\frac{\Gamma_q(\lambda+1)}{\Gamma_q(\lambda+1+\alpha)}\
x^{\lambda+\alpha}(a/x;q)_{\lambda+\alpha}\qquad (a<x)\ ,\\
D_{q,a}^{\alpha}\bigl(x^{\lambda}(a/x;q)_{\lambda}\bigr) &=&
\frac{\Gamma_q(\lambda+1)}{\Gamma_q(\lambda+1-\alpha)}\
x^{\lambda-\alpha}(a/x;q)_{\lambda-\alpha} \ , \\ \\
{}_\star D_{q,a}^{\alpha}\bigl(x^{\lambda}(a/x;q)_{\lambda}\bigr) &=&
\left\{
\begin{array}{cc} 0\ , & \lambda\in \mathbb N_0; \ \alpha>\lambda \ ,\\    \\
D_{q,a}^{\alpha}\bigl(x^{\lambda}(a/x;q)_{\lambda}\bigr) ,&
\text{\rm otherwise.}
\end{array}
\right.
\end{eqnarray*}
\end{lem}

\noindent{\it Proof.} For $\lambda\ne0$, according to the
definition (\ref{fracAgain}), we have
$$
\aligned
& I_{q,a}^{\alpha}\bigl(x^{\lambda}(a/x;q)_{\lambda}\bigr)\\
&\qquad =\frac{x^{\alpha-1}}{\Gamma_q(\alpha)} \Bigl(\int_0^x
(qt/x;q)_{\alpha-1} t^{\lambda}(a/t;q)_{\lambda} d_qt -\int_0^a
(qt/x;q)_{\alpha-1} t^{\lambda}(a/t;q)_{\lambda} d_qt\Bigr).
\endaligned
$$
Also, the following is valid:
$$
\int_0^a (qt/x;q)_{\alpha-1}\ t^{\lambda}(a/t;q)_{\lambda} d_qt
=a^{\lambda+1}(1-q)\sum_{k=0}^\infty (aq^{k+1}/x;q)_{\alpha-1}\
q^{k\lambda}(q^{-k};q)_{\lambda}\ q^k\ ,
$$
what vanishes because of (\ref{fracprop4}). Therefrom, according
to definition (\ref{(qint0)}), we get
$$
\aligned
\int_0^x (qt/x;q)_{\alpha-1}\ &t^{\lambda}(a/t;q)_{\lambda}\ d_qt\\
&=x^{\lambda+1}(1-q)\sum_{k=0}^\infty (q^{1+k};q)_{\alpha-1}\
\bigl(a/(xq^k);q\bigr)_{\lambda}\ q^{(\lambda+1)k}\ .
\endaligned
$$
We notice presence of (\ref{distrib1}) in the previous formula,
i.e.
$$
\aligned
\int_0^x (qt/x;q)_{\alpha-1}\ &t^{\lambda}(a/t;q)_{\lambda}\ d_qt\\
&=(1-q)\ x^{\lambda+1}\ (q;q)_{\alpha-1}\ (q;q)_{\lambda}\ S\Bigl(
\lambda+1,\alpha,a/(qx)\Bigr)\ .
\endaligned
$$
By using (\ref{distrib}), we get
$$
\int_0^x (qt/x;q)_{\alpha-1} t^{\lambda}(a/t;q)_{\lambda} d_qt
=(1-q)\frac{(q;q)_{\alpha-1}(q;q)_{\lambda}}{(q;q)_{\alpha+\lambda}}\
x^{\lambda+1}(a/x;q)_{\alpha+\lambda},
$$
and applying (\ref{Gammaq1}), we obtain the required formula for
$I_{q,a}^{\alpha}\bigl(x^{\lambda}(a/x;q)_{\lambda}\bigr)$ when
$\lambda\ne0$.

In the case when $\lambda = 0$, using $q$--integration by parts
(\ref{IntParts}), we have
$$
\aligned (I_{q,a}^{\alpha}{\bf
1})(x)&=\frac{x^{\alpha-1}}{\Gamma_q(\alpha)}\int_a^x
(qt/x;q)_{\alpha-1}d_qt =\frac{1}{\Gamma_q(\alpha)}
\int_a^x \frac{D_q\bigl(x^{\alpha}(t/x;q)_{\alpha}\bigr)}{-[\alpha]_q}d_qt\\
&=\frac{-1}{\Gamma_q(\alpha+1)}\int_a^x
D_q\bigl(x^{\alpha}(t/x;q)_{\alpha}\bigr)d_qt
=\frac{1}{\Gamma_q(\alpha+1)}x^{\alpha}(a/x;q)_{\alpha} \ .
\endaligned
$$
The terms for $q$--derivatives can be obtained by applying
definitions (\ref{Rie-Liou}) and (\ref{Caputo}).\ \ $\Box$

\begin{cor}\label{Irec8}
For $\alpha\in\mathbb R^{+}\setminus \mathbb N_0$, $n\in\mathbb
N_0$, and $a<x$, the following is valid:
\begin{eqnarray*}
I_{q,a}^{\alpha}\bigl(x^{n}\bigr)&=&(1-q)^\alpha \sum_{k=0}^n
a^{n-k}(q^{n-k+1};q)_k\ \frac{x^{k+\alpha}(a/x;q)_{k+\alpha}}{(q;q)_{k+\alpha}}\ . \\
D_{q,a}^{\alpha}(x^n) &=&(1-q)^{-\alpha} \sum_{k=0}^n
a^{n-k}(q^{n-k+1};q)_k\ \frac{x^{k-\alpha}(a/x;q)_{k-\alpha}}{(q;q)_{k-\alpha}}\ . \\
 {}_\star D_{q,a}^{\alpha}(x^n)
 &=&\frac{(q^{n+1-\lceil\alpha\rceil};q)_{\lceil\alpha\rceil}}
{(1-q)^{\alpha}}\sum_{k=\lceil\alpha\rceil}^{n}
a^{n-k}(q^{n-k+1};q)_k\
\frac{x^{k-\alpha}(a/x;q)_{k-\alpha}}{(q;q)_{k-\alpha}}\ .
\end{eqnarray*}

(Notice that ${}_\star D_{q,a}^{\alpha}(x^n)=0$ when $\alpha>n$.)
\end{cor}

The $q$--exponential functions (see \cite{BHS}) can be written
like power series or, applying $q$--form of Taylor theorem
(\ref{qTaylor}), by
\begin{eqnarray}
e_q(x)&=&\sum_{n=0}^{\infty}\frac{x^n}{(q;q)_n}  =
e_q(a)\sum_{n=0}^\infty \frac{x^{n}(a/x;q)_{n}}{(q;q)_{n}} \
\qquad \ (|x|<1)\ ,\label{eqx}\\
E_q(x)&=&\sum_{n=0}^{\infty}\frac{q^{\binom n2}}{(q;q)_n}\
x^n=E_q(a)\sum_{n=0}^\infty \frac{q^{\binom n2}}{(-a;q)_n}\
\frac{x^{n}(a/x;q)_{n}}{(q;q)_{n}}\ .\label{vEqx}
\end{eqnarray}

\begin{cor}\label{Ieqx}
For $\alpha\in\mathbb R^{+}\setminus \mathbb N_0$ and $0<a<x<1$,
the following is valid:
\begin{eqnarray*}
I_{q,a}^{\alpha}\bigl(e_q(x)\bigr) &=&(1-q)^{\alpha}\
e_q(a)\sum_{n=0}^\infty
\frac{x^{n+\alpha}(a/x;q)_{n+\alpha}}{(q;q)_{n+\alpha}}\ ,\\
D_{q,a}^{\alpha}(e_q(x))& =& (1-q)^{-\alpha}\
e_q(a)\sum_{n=0}^\infty
\frac{x^{n-\alpha}(a/x;q)_{n-\alpha}}{(q;q)_{n-\alpha}}\ , \\
{}_\star D_{q,a}^{\alpha}(e_q(x))& =& (1-q)^{-\alpha}\
e_q(a)\sum_{n=\lceil\alpha\rceil}^\infty
\frac{x^{n-\alpha}(a/x;q)_{n-\alpha}}{(q;q)_{n-\alpha}}\ .
\end{eqnarray*}
\end{cor}

\begin{cor}\label{IEqx}
For $\alpha\in\mathbb R^{+}\setminus \mathbb N_0$ and $0<a<x$, the
following is valid:
\begin{eqnarray*}
I_{q,a}^{\alpha}\bigl(E_q(x)\bigr) &=&(1-q)^{\alpha}\
E_q(a)\sum_{n=0}^\infty \frac{q^{\binom n2}}{(-a;q)_n}\
\frac{x^{n+\alpha}(a/x;q)_{n+\alpha}}{(q;q)_{n+\alpha}}\ ,\\
D_{q,a}^{\alpha}\bigl(E_q(x)\bigr)
 &=&\frac{E_q(a)}{(1-q)^{\alpha}}\sum_{n=0}^\infty \frac{q^{\binom n2}}{(-a;q)_n}\
\frac{x^{n-\alpha}(a/x;q)_{n-\alpha}}{(q;q)_{n-\alpha}}\ , \\
{}_\star D_{q,a}^{\alpha}\bigl(E_q(x)\bigr)
 &=&\frac{E_q(a)}{(1-q)^{\alpha}}\sum_{n=\lceil\alpha\rceil}^\infty \frac{q^{\binom n2}}{(-a;q)_n}\
\frac{x^{n-\alpha}(a/x;q)_{n-\alpha}}{(q;q)_{n-\alpha}}\ .
\end{eqnarray*}
\end{cor}

\section{The relationship between fractional \\ $q$--integrals and $q$--derivatives}

It is very important to establish the connection between two types
of the fractional $q$--derivatives.

\begin{thm}\label{RL-C}
Let $\alpha\in \mathbb R^+\setminus\mathbb N_0$ and $a<x$. The
connection between Caputo type and Riemann-Liouville type
fractional integral is
$$
\bigl( D_{q,a}^{\alpha} f\bigr)(x)
=\bigl({}_{\star}D_{q,a}^{\alpha}
f\bigr)(x)+\sum_{k=0}^{\lceil\alpha\rceil-1}\frac{(D_q^k
f)(a)}{\Gamma_q(1+k-\alpha)}x^{k-\alpha}(a/x;q)_{k-\alpha}
$$
\end{thm}
\noindent{\it Proof.} Any $\alpha\in \mathbb R^+\setminus\mathbb
N_0$ we can write in the form $\alpha=n+\varepsilon$, where
$\varepsilon\in (0,\ 1)$. We will prove the statement by
mathematical induction over $n\in\mathbb N_0$.

At first, let $n=0$, i.e., $\alpha\in(0,1)$. According to Lemma
\ref{Ialbasic}, we have
$$
\aligned
 \bigl( I_{q,a}^{1-\alpha}f\bigr)(x)
&=\bigl(I_{q,a}^{2-\alpha} D_q
f\bigr)(x)+\frac{f(a)}{\Gamma_q(2-\alpha)}
x^{1-\alpha}(a/x;q)_{1-\alpha}\\
&=\Bigl(I_{q,a} \bigl({}_{\star}D^{\alpha}_{q,a}
f\bigr)\Bigr)(x)+\frac{f(a)}{\Gamma_q(2-\alpha)}x^{1-\alpha}(a/x;q)_{1-\alpha}\
.
\endaligned
$$
By $q$--deriving, we get
$$
 \bigl( D_qI_{q,a}^{1-\alpha}f\bigr)(x)=
\Bigl(D_qI_{q,a} \bigl({}_{\star}D^{\alpha}_{q,a} f\bigr)\Bigr)(x)
+
\frac{f(a)}{\Gamma_q(2-\alpha)}D_q\Bigl(x^{1-\alpha}(a/x;q)_{1-\alpha}\Bigr),
$$
and, with respect to (\ref{xminusa}),
$$
 \bigl( D_q^{\alpha}f\bigr)(x)=
\bigl({}_{\star}D^{\alpha}_{q,a} f\bigr)(x) +
\frac{f(a)}{\Gamma_q(1-\alpha)}x^{-\alpha}(a/x;q)_{-\alpha}\ .
$$
Suppose that the statement is valid for a real
$\alpha=n+\varepsilon$, $\varepsilon\in (0,1)$, for a positive
integer $n\in\mathbb N$ and let us prove that it is valid for
$\alpha=n+1+\varepsilon$. Indeed, according to Theorem
\ref{DRalmE}, the next equality is valid:
$$
\bigl(D^{\alpha}_{q,a}f\bigr)(x)=\Bigl( D_q
D_{q,a}^{n+\varepsilon}f\Bigr)(x).
$$
With respect to the inductional assumption
$$
\bigl( D_{q,a}^{n+\varepsilon} f\bigr)(x)
=\bigl({}_{\star}D_{q,a}^{n+\varepsilon}
f\bigr)(x)+\sum_{k=0}^{n}\frac{(D_q^k
f)(a)}{\Gamma_q(1+k-n-\varepsilon)}x^{k-n-\varepsilon}(a/x;q)_{k-n-\varepsilon},
$$
and the formula (\ref{xminusa}),  we can write
$$
\aligned \bigl(D^{\alpha}_{q,a}f\bigr)&(x) \\
&= \bigl(D_q \
{}_{\star}D_{q,a}^{n+\varepsilon}f\bigr)(x)
+\sum_{k=0}^{n}\frac{(D_q^k f)(a)}{\Gamma_q(1+k-n-\varepsilon)}
D_q\bigl(x^{k-n-\varepsilon}(a/x;q)_{k-n-\varepsilon}\bigr)\\
&= \bigl(D_q \ {}_{\star}D_{q,a}^{n+\varepsilon} f\bigr)(x)
+\sum_{k=0}^{n}\frac{(D_q^k f)(a)}{\Gamma_q(k-n-\varepsilon)} \
x^{k-n-1-\varepsilon}(a/x;q)_{k-n-1-\varepsilon} \ .
\endaligned
$$
Using the Theorem \ref{DCalmD}, we obtain
$$
\aligned \bigl(D_q \ {}_{\star}D_{q,a}^{n+\varepsilon} f\bigr)(x)
=\bigl({}_{\star}D_{q,a}^{n+1+\varepsilon} f\bigr)(x)+
\dfrac{\bigl(D_{q}^{n+1}f\bigr)(a)}{\Gamma_q(1-\varepsilon)}\
x^{-\varepsilon}(a/x;q)_{-\varepsilon}\ .
\endaligned
$$
So,
$$
\aligned \bigl(D^{\alpha}_{q,a}f\bigr)(x)
&=\bigl({}_{\star}D_{q,a}^{n+1+\varepsilon} f\bigr)(x)+
\dfrac{\bigl(D_{q}^{n+1}f\bigr)(a)}{\Gamma_q(1-\varepsilon)}\
x^{-\varepsilon}(a/x;q)_{-\varepsilon}\\
&+\sum_{k=0}^{n}\frac{(D_q^k f)(a)}{\Gamma_q(k-n-\varepsilon)}
\ x^{k-n-1-\varepsilon}(a/x;q)_{k-n-1-\varepsilon}\\
&=\bigl({}_{\star}D_{q,a}^{\alpha} f\bigr)(x)+
\sum_{k=0}^{n+1}\frac{(D_q^k f)(a)}{\Gamma_q(k-n-\varepsilon)} \
x^{k-n-1-\varepsilon}(a/x;q)_{k-n-1-\varepsilon}\ ,
\endaligned
$$
what is finishing the proof.\ \ $\Box$

\smallskip

Here, we will discuss behavior of compositions of previously
defined operators.

\begin{thm} \label{9}
Let $\alpha\in \mathbb R^+$. Then, for $a<x$, the following is
valid:
$$
\bigl(D_{q,a}^{\alpha}I_{q,a}^{\alpha} f\bigr)(x) = f(x)\ .
$$
\end{thm}
\noindent{\it Proof.} With respect to Theorem \ref{kom} and the
formulas (\ref{DerInt}) and (\ref{IntDer}), we have
$$
\aligned
\bigl(D_{q,a}^{\alpha}I_{q,a}^{\alpha} f\bigr)(x)
&=\bigl(D^{\lceil\alpha\rceil}_qI_{q,a}^{\lceil\alpha\rceil-\alpha}I_{q,a}^{\alpha}f\bigr)(x)
=\bigl(D^{\lceil\alpha\rceil}_qI_{q,a}^{\lceil\alpha\rceil-\alpha+\alpha}f\bigr)(x)\\
&=\bigl(D^{\lceil\alpha\rceil}_qI_{q,a}^{\lceil\alpha\rceil}f\bigr)(x)=f(x)\
.\ \ \Box
\endaligned
$$

\begin{thm}\label{K}
Let $\alpha\in\mathbb R^+\setminus\mathbb N$. Then
$$
\bigl(I_{q,a}^{\alpha}D_{q,a}^{\alpha} f\bigr)(x) = f(x)\qquad
(a<x)\ .
$$
\end{thm}

\noindent{\it Proof.} Let $\alpha\in(0,1)$. Since, according to
(\ref{IntDer}), we can write
$$
f(x)=\bigl(I_{q,a}D_q f\bigr)(x) + f(a)\ ,
$$
and, by using Theorem \ref{kom} and Lemma \ref{Irec} we have
$$
\aligned \bigl(I_{q,a}^{1-\alpha}f\bigr)(x) &=
\bigl(I_{q,a}^{1-\alpha}I_{q,a}D_q f\bigr)(x)
+f(a)\bigl(I_{q,a}^{1-\alpha} {\bf 1}\bigr)(x)\\
& =\bigl(I_{q,a}^{2-\alpha}D_q f\bigr)(x)
+\frac{f(a)}{\Gamma_q(2-\alpha)}x^{1-\alpha}(a/x;q)_{1-\alpha}\ .
\endaligned
$$
Applying $D_q$ on both sides of equality, we obtain
$$
\aligned \bigl(D_{q,a}^{\alpha}f\bigr)(x)&=\bigl(D_q
I_{q,a}^{1-\alpha}f\bigr)(x)\\
& =\bigl(D_q I_{q,a}^{2-\alpha}D_q f\bigr)(x)
+\frac{f(a)}{\Gamma_q(2-\alpha)}D_q\bigl(x^{1-\alpha}(a/x;q)_{1-\alpha}\bigr)\\
&=\bigl(I_{q,a}^{1-\alpha}D_q f\bigr)(x)
+\frac{f(a)}{\Gamma_q(1-\alpha)}x^{-\alpha}(a/x;q)_{-\alpha}\ .
\endaligned
$$
Now, again with respect to Theorem \ref{kom} and Lemma \ref{Irec},
the following is valid:
$$
\aligned \bigl(I_{q,a}^{\alpha}D_{q,a}^{\alpha}f\bigr)(x)
&=\bigl(I_{q,a}^{\alpha} I_{q,a}^{1-\alpha}D_q f\bigr)(x)
+\frac{f(a)}{\Gamma_q(1-\alpha)}I_{q,a}^{\alpha}\bigl(x^{-\alpha}(a/x;q)_{-\alpha})\\
&=\bigl(I_{q,a}D_q f\bigr)(x)+f(a)=f(x)\ .
\endaligned
$$
Let $\alpha=n+\varepsilon$, with $n\in\mathbb N$, $0<
\varepsilon<1$. Putting $\alpha\mapsto \alpha-1$ and $f\mapsto
D_{q,a}^{\alpha-1}f$ into Lemma \ref{Ialbasic}, and applying
Theorem~\ref{DRalmE}, we get
$$
\aligned (I_{q,a}^{\alpha-1}D_{q,a}^{\alpha-1}f)(x)
&=(I_{q,a}^{\alpha}D_q D_{q,a}^{\alpha-1}f)(x)
+\frac{(D_{q,a}^{\alpha-1}f)(a)}{\Gamma_q(\alpha)}x^{\alpha-1}(a/x;q)_{\alpha-1}\\
&=(I_{q,a}^{\alpha}D_{q,a}^{\alpha}f)(x)
+\frac{(D_{q,a}^{\alpha-1}f)(a)}{\Gamma_q(\alpha)}x^{\alpha-1}(a/x;q)_{\alpha-1}\
.
\endaligned
$$
According to property (\ref{Dqf(a)}), we conclude that
$$
(I_{q,a}^{\alpha}D_{q,a}^{\alpha}f)(x)=
(I_{q,a}^{\alpha-1}D_{q,a}^{\alpha-1}f)(x)\ .
$$
Repeating the last identity $n$ times, we get
$$
(I_{q,a}^{\alpha}D_{q,a}^{\alpha}f)(x)=
(I_{q,a}^{\alpha-n}D_{q,a}^{\alpha-n}f)(x)=
(I_{q,a}^{\varepsilon}D_{q,a}^{\varepsilon}f)(x)=f(x)\ ,
$$
what is finishing the proof.\ \ $\Box$

\begin{thm} \label{17}
Let $\alpha\in \mathbb R^+\setminus\mathbb N$. Then, for $a<x$,
the following is valid:
$$
\bigl(I_{q,a}^{\alpha}\ {}_{\star}D_{q,a}^{\alpha} f\bigr)(x) =
f(x)-\sum_{k=0}^{\lceil\alpha\rceil -1}\dfrac{\bigl(D_q^k
f\bigr)(a)}{[k]_q!}\ x^k(a/x;q)_k\ .
$$
\end{thm}

\noindent{\it Proof.} With respect to Theorem \ref{kom} and the
formulas (\ref{DerInt}) and (\ref{IntDer}), we have
$$
\aligned \bigl(I_{q,a}^{\alpha}\ {}_{\star}D_{q,a}^{\alpha}
f\bigr)(x)
&=\bigl(I_{q,a}^{\alpha}I_{q,a}^{\lceil\alpha\rceil-\alpha}D^{\lceil\alpha\rceil}_q
f\bigr)(x)
=\bigl(I_{q,a}^{\lceil\alpha\rceil}D^{\lceil\alpha\rceil}_qf\bigr)(x)\\
&=f(x)-\sum_{k=0}^{\lceil\alpha\rceil -1}\dfrac{\bigl(D_q^k
f\bigr)(a)}{[k]_q!}\ x^k(a/x;q)_k\ .\ \Box
\endaligned
$$

\begin{thm} \label{18}
Let $\alpha\in \mathbb R^+\setminus\mathbb N$. Then, for $a<x$,
the following is valid:
$$
\bigl({}_{\star}D_{q,a}^{\alpha}I_{q,a}^{\alpha} f\bigr)(x) =
f(x)\ .
$$
\end{thm}

\noindent{\it Proof.} Putting $f\mapsto I_{q,a}^{\alpha}f$ into
Theorem \ref{RL-C}, and using Theorem~\ref{K},
Corollary~\ref{DIal} and formula (\ref{qFInt0}), we get
$$
\aligned
\bigl({}_{\star}D_{q,a}^{\alpha}I_{q,a}^{\alpha} f\bigr)(x)
&=\bigl(D_{q,a}^{\alpha} I_{q,a}^{\alpha} f\bigr)(x)-\sum_{k=0}^{\lceil\alpha\rceil-1}
\frac{(D_q^k I_{q,a}^{\alpha}f)(a)}{\Gamma_q(1+k-\alpha)}x^{k-\alpha}(a/x;q)_{k-\alpha}\\
&=f(x)-\sum_{k=0}^{\lceil\alpha\rceil-1}
\frac{(I_{q,a}^{\alpha-k}f)(a)}{\Gamma_q(1+k-\alpha)}x^{k-\alpha}(a/x;q)_{k-\alpha}=f(x)\ .
\ \Box
\endaligned
$$

\begin{thm} \label{DqIqAlBe}
Let $\alpha \in \mathbb R$ and $\beta \in \mathbb R^+$. Then, for
$a<x$, the following is valid:
$$
\bigl(D_{q,a}^{\alpha} I_{q,a}^{\beta} f\bigr)(x) =
\bigl(D_{q,a}^{\alpha-\beta} f\bigr)(x)\ .
$$
\end{thm}

\noindent{\it Proof.} Let $\alpha = n+\varepsilon$ and $\beta =
m+\delta$, where $n>m$ and  $\varepsilon, \delta\in [0,1)$ such
$\varepsilon< \delta$. Then
$$
\aligned \bigl(D_{q,a}^{\alpha} \ I_{q,a}^{\beta} f\bigr)(x)
&=\bigl(D_{q,a}^{n+1} \ I_{q,a}^{1-\varepsilon} \
I_{q,a}^{m+\delta} f\bigr)(x)
\\
&=\bigl(D_{q}^{n+1} \ I_{q,a}^{m+1+\delta-\varepsilon} \
 f\bigr)(x) \\
 &=\bigl(D_{q}^{n+1} \ I_{q,a}^{m+1} \ I_{q,a}^{\delta-\varepsilon} \
 f\bigr)(x) \\
 &=\bigl(D_{q}^{n-m} \ I_{q,a}^{\delta-\varepsilon}
 f\bigr)(x)  \ .
 \endaligned
$$
From the other side
$$
\bigl(D_{q,a}^{\alpha-\beta} f\bigr)(x) =
\bigl(D_{q}^{\lceil\alpha-\beta\rceil} \
I_{q,a}^{\lceil\alpha-\beta\rceil-(\alpha-\beta)} f\bigr)(x)
 = \bigl(D_{q}^{n-m} \
I_{q,a}^{\delta-\varepsilon} f\bigr)(x) \ . \ \Box
$$
\begin{thm} \label{IqDqAlBe}
Let $\alpha \in \mathbb R\setminus\mathbb N$ and $\beta \in
\mathbb R^+$. Then, for $a<x$, the following is valid:
$$
\bigl(I_{q,a}^{\beta} D_{q,a}^{\alpha} f\bigr)(x) =
\bigl(D_{q,a}^{\alpha-\beta} f\bigr)(x)\ .
$$
\end{thm}

\noindent{\it Proof.} If $\alpha\le 0$, the statement follows immediately from
the definition and Theorem~\ref{kom}.

Let $0<\alpha\le \beta$. Then, with respect to Theorem~\ref{kom} and
Theorem~\ref{K}, we have
$$
\bigl(I_{q,a}^{\beta} D_{q,a}^{\alpha} f\bigr)(x)
=\bigl(I_{q,a}^{\beta-\alpha}I_{q,a}^{\alpha} D_{q,a}^{\alpha}
f\bigr)(x)
=\bigl(I_{q,a}^{\beta-\alpha} f\bigr)(x)
=\bigl(D_{q,a}^{\alpha-\beta} f\bigr)(x)\ .
$$

Finally, let $\alpha>\beta$. According to Theorem~\ref{K}, we can write
$$
f(x)=\bigl(I_{q,a}^{\alpha} D_{q,a}^{\alpha} f\bigr)(x)=
\bigl(I_{q,a}^{\alpha-\beta}I_{q,a}^{\beta} D_{q,a}^{\alpha} f\bigr)(x) \ .
$$
Applying $D_{q,a}^{\alpha-\beta}$ on both sides of the last equality,
we finish the proof.\ \ $\Box$

\medskip

Notice that statement of Theorem~\ref{IqDqAlBe} is not valid for
$\alpha\in\mathbb N$. In that case, the following identity holds:
$$
\bigl(I_{q,a}^{\beta} D_{q}^{n} f\bigr)(x) =
\bigl(D_{q,a}^{n-\beta} f\bigr)(x)
-\sum_{k=0}^{n-1}\frac{\bigl(D_q^k f\bigr)(a)}{\Gamma_q(\beta-n+k+1)}
\ x^{\beta-n+k}(a/x;q)_{\beta-n+k}\ .
$$
Indeed, if $\alpha=n\le \beta$, by using Theorem~\ref{kom},
formula (\ref{IntDer}) and Corollary~\ref{Irec}, we get
$$
\aligned
\bigl(I_{q,a}^{\beta} D_{q}^{n} f\bigr)(x)
&=\bigl(I_{q,a}^{\beta-n}I_{q,a}^{n} D_{q}^{n} f\bigr)(x)\\
&=\bigl(I_{q,a}^{\beta-n} f\bigr)(x)-
\sum_{k=0}^{n-1}\frac{\bigl(D_q^k f\bigr)(a)}{[k]_q!}
I_{q,a}^{\beta-n}\bigl( x^{k}(a/x;q)_{k}\bigr)\\
&=\bigl(D_{q,a}^{n-\beta} f\bigr)(x)-
\sum_{k=0}^{n-1}\frac{\bigl(D_q^k f\bigr)(a)}{\Gamma_q(\beta-n+k+1)}
x^{\beta-n+k}(a/x;q)_{\beta-n+k}\ .
\endaligned
$$

In similar way, by using Theorem~\ref{RL-C}, Theorem~\ref{K},
Theorem~\ref{DqIqAlBe} and Theorem~\ref{IqDqAlBe}, the next
properties can be proven.

\begin{thm} \label{*DqIqAlBe}
Let $\alpha \in \mathbb R\setminus\mathbb N$ and $\beta \in
\mathbb R^+$. Then, for $a<x$, the following is valid:
$$
\aligned
\bigl({}_\star D_{q,a}^{\alpha} & I_{q,a}^{\beta} f\bigr)(x)\\
&=\bigl({}_\star D_{q,a}^{\alpha-\beta} f\bigr)(x)
+\sum_{k=0}^{\lceil\alpha-\beta\rceil-1}
\frac{\bigl(D_q^k f\bigr)(a)}{\Gamma_q(k-\alpha+\beta+1)}
x^{k-\alpha+\beta}(a/x;q)_{k-\alpha+\beta}\ ,\\
\bigl(I_{q,a}^{\beta}\ &{}_\star D_{q,a}^{\alpha}f\bigr)(x)\\
&=\bigl({}_\star D_{q,a}^{\alpha-\beta} f\bigr)(x)
-\sum_{k=\lceil\alpha-\beta\rceil}^{\lceil\alpha\rceil-1}
\frac{\bigl(D_q^k f\bigr)(a)}{\Gamma_q(k-\alpha+\beta+1)}
x^{k-\alpha+\beta}(a/x;q)_{k-\alpha+\beta}\ .
\endaligned
$$
\end{thm}

\begin{thm} \label{IqDq}
Let $a\le c<x$ and $\alpha\in \mathbb R^+\setminus\mathbb N$. Then
the following is valid:
$$
\aligned \bigl(I_{q,c}^{\alpha}D_{q,a}^{\alpha} f\bigr)(x) &=
\bigl(I_{q,c}^{\alpha-\lceil\alpha\rceil+1}D_{q,a}^{\alpha-\lceil\alpha\rceil+1}
f\bigr)(x)\\
&-\sum_{k=1}^{\lceil\alpha\rceil-1}
\frac{\bigl(D_{q,a}^{\alpha-k}f\bigr)(c)}{\Gamma_q(\alpha-k+1)}\
x^{\alpha-k}(c/x;q)_{\alpha-k} \ .
\endaligned
$$
\end{thm}

\medskip

\noindent{\bf Acknowledgements}

\smallskip

 This work was supported by Ministry of Science, Technology and Development of
Republic Serbia, through the project No $144023$ and No $144013$.

\end{document}